\documentclass[a4paper]{article}

\usepackage[english]{babel}
\usepackage[utf8]{inputenc}
\usepackage{amsmath}
\usepackage{graphicx}
\usepackage{amssymb}
\usepackage{amsthm}
\usepackage{tikz-cd}
\usepackage{mathrsfs}
\usepackage[colorinlistoftodos]{todonotes}
\usepackage{enumitem}
\usepackage{yfonts}
\usepackage{mathtools}
\usepackage{hyperref}

\usepackage[isbn=false,url=false,eprint=false]{biblatex}

\title{Subconvexity for Symmetric Square Off-Centre}
\author{Mayukh Dasaratharaman, Ritabrata Munshi}

\newtheorem{thm}{Theorem}
\newtheorem{lem}{Lemma}

\newtheorem*{thm*}{Theorem}
\newtheorem*{lem*}{Lemma}
\newtheorem*{quest*}{Question}
\newtheorem*{defn*}{Definition}
\newtheorem*{eg*}{Example}
\newtheorem*{ex*}{Exercise}
\newtheorem*{conj*}{Conjecture}
\newtheorem*{cor*}{Corollary}
\newtheorem*{rmk*}{Remark}
\newtheorem*{prop*}{Proposition}

\newcommand{\R}{\mathbb{R}}

\newcommand{\Z}{\mathbb{Z}}

\begin{document}
	\maketitle
	\begin{abstract}
		Let $p$ be a prime. Let $f$ be a holomorphic modular form of level $p$ with trivial nebentypus. We prove the bound
		
		$$L\left(\textup{sym}^2f, \frac{1}{2} + it\right) \ll_{f,\epsilon} p^{1/2+\epsilon}t^{3/4-1/12 + \epsilon}$$
		
		This bound is subconvex in the $t$-aspect and almost convex in the level aspect simultaneously.
	\end{abstract}
	\section{Introduction}
	Let $p$ be a prime. Let $f$ be a holomorphic new form of level $p$ and trivial nebentypus. Consider the symmetric square L-function associated to $f$ defined by
	
	$$L(s,\textup{sym}^2 f) = \zeta^{(p)}(2s)\sum \limits_{n =1}^\infty \frac{a_f(n^2)}{n^s}$$
	
	valid for $\Re(s) > 1$. Here $\zeta^{(p)}(2s)$ is the Riemann zeta function with the $p$-th Euler factor missing and $\lambda_f(n)$ are the normalised Fourier coefficients associated to $f$. The L-function also has an Euler product $L(s,\textup{sym}^2 f) = \prod \limits_q L_q(s, \textup{sym}^2 f)$ where
	
	$$L_q(s,\textup{sym}^2 f) = \left(1 - \frac{\lambda_f(q^2)}{q^s} + \frac{\lambda_f(q^2)}{q^{2s}} - \frac{1}{p^{3s}} \right)^{-1}$$
	
	if $q \neq p$ and
	
	$$L_p(s, \textup{sym}^2 f) = \left(1 - \frac{1}{p^{s+1}} \right)^{-1}$$
	
	We define the completed L-function
	
	$$\Lambda(s,\textup{sym}^2 f) = p^{s/2} L_\infty(s)L(s,\textup{sym}^2 f)$$
	
	where
	
	$$L_\infty(s) = \pi^{- \frac{3s}{2}} \Gamma \left( \frac{s+1}{2} \right) \Gamma \left( \frac{s+k-1}{2} \right) \Gamma \left(\frac{s+k}{2} \right)$$
	
	Then $\Lambda(s,\textup{sym}^2 f)$ is an entire function and satisfies the functional equation $\Lambda(s,\textup{sym}^2 f) = \Lambda(1-s,\textup{sym}^2 f)$.
	
	An important problem is to study the growth of $L(s,f)$ on the critical line $\Re(s) = 1/2$, say with $s=1/2+it$ in terms of $t$ and $p$. From the functional equation we can prove that
	
	$$L\left(\frac{1}{2} + it,\textup{sym}^2f \right) \ll p^{1/2+\epsilon}t^{3/4+\epsilon}$$
	
	which is known as the convexity estimate while the Lindelof hypothesis predicts that
	
	$$L\left(\frac{1}{2} + it,\textup{sym}^2f \right) \ll_{p,t} (pt)^\epsilon$$
	
	for any $\epsilon > 0$. Proving any bound that is better than the convexity bound unconditionally is a difficult problem.
	We prove the following hybrid subconvexity bound given in the following theorem.
	
	\begin{thm}\label{thm1}
		
		We have for any $\epsilon > 0$
		
		$$L(\frac{1}{2} + it, \textup{sym}^2f) \ll_{\epsilon} p^{1/2+\epsilon}t^{3/4 - 1/12 + \epsilon}$$
		
	\end{thm}
	
	The theorem states that we have a bound that is almost convex in the level aspect and subconvex in the $t$-aspect simultaneously. This directly improves the bound of Iwaniec-Michel if we consider only the level and $t$-aspects simultaneously. However, we note that the exponent in the $t$-aspect is not best possible. In \cite{lin_strong_2022} the authors obtain

    $$L \left( \frac{1}{2}+it, \textup{sym}^2 f \right) \ll t^{3/4-3/20 + \epsilon}$$

    See also \cite{pal} which makes improvements.
	
	The symmetric square L-function is an example of a degree 3 L-function. The first known subconvex result for degree 3 was proved by Li \cite{li_bounds_2011} in the $t$-aspect. Building on the work of Conrey and Iwaniec \cite{conrey_iwaniec}, she obtained the bound
	
	$$L(\frac{1}{2} + it,\pi) \ll t^{3/4 - 1/16 + \epsilon}$$
	
	where $\pi$ is a symmetric square or a quadratic twist of a symmetric square. The proof required the form to be self-dual, and the techniques couldn't be extended to a generic degree 3 L-function.(see also \cite{lin_strong_2022}, \cite{nunes_subconvexity_2017},  for subconvex bounds with similar proof). Munshi later proved a subconvex bound
	
	$$L \left( \frac{1}{2} + it,\pi \right) \ll t^{3/4-1/16 + \epsilon}$$
	
	The proof used the delta method which didn't require the L-function to be self-dual and could be extended to a generic degree 3 form. The exponent in the bound has recently been improved to $3/4-1/8 + \epsilon$ by Aggarwal, Leung and Munshi \cite{aggarwal_short_2022}.
	
	For a degree $n$ automorphic form $\pi_n$ in the generic position, Nelson \cite{nelson_bounds_2023} proved that
	
	$$L\left( \frac{1}{2} + it,\pi_n \right) \ll t^{(1-\delta_n)n/4}$$
	
	with $\delta_n \geq \frac{1}{3n^5 - 2n^4 - n^2}$. Note that though the bounds are weaker than the current best known for $n \leq 3$, it holds for any $n$.
	
	\section{Tools and Outline}
	
	By the approximate functional equation and a smooth partition of unity, we have (see \cite{iwaniec_second_2001})
	
	$$L\left(\frac{1}{2} +it,\textup{sym}^2 f \right)\ll (pt)^\epsilon \sup \limits_{N \ll p^{1+\epsilon}t^{3/2+\epsilon}} \frac{S(N)}{\sqrt{N}}$$
	
	where
	
	\begin{equation}\label{defnS(N)}
		S(N) = \sum \limits_{n} \lambda_f(n^2)n^{-it}V\left( \frac{n}{N}\right)
	\end{equation}
	
	where $V(x)$ is a smooth function with compact support in $[1/2,5/2]$, is identically 1 on $[1,2]$ and satisfies the bounds $V^{(j)}(x) \ll_j 1$.
	
	We now wish to separate the oscillation of $\lambda_f(n^2)$ from $m^{-it}$. To do this we apply the following delta symbol (see \cite{munshi_subconvexity_2021})
	
	\begin{equation}
		\delta(m=n^2) = \frac{1}{Q} \sum \limits_{q \leq Q} \int \limits_\R g(q,x) \sideset{}{^*} \sum \limits_{ a \bmod q} \frac{1}{q} e \left( \frac{a(n^2-m)}{q}\right) e \left(\frac{(n^2-m)x}{qQ}\right) dx
	\end{equation}
	
	valid for any $Q \geq 1$. Here $g(q,x)$ satisfies
	
	\begin{equation}{\label{smallq}}
		g(q,x) = 1 + O\left(\frac{1}{qQ} \left( \frac{q}{Q} + |x|\right)^A\right)\end{equation}
	\begin{equation}{\label{smallx}} g(q,x) \ll_A |x|^{-A}\end{equation}
	
	for any $A > 0$.and for every $k \geq 1$,

	Thus, $(\ref{smallx})$ says for all $q$, the integral is essentially supported on $[-Q^\epsilon,Q^\epsilon]$ and $(\ref{smallq})$ says can replace $g(q,x)$ by $1$ at the cost of a negligible error term for $q \ll Q^{1-\epsilon}$. For $q \gg Q^{1-\epsilon}$, while we don't know that $g(q,x)$ is $1$, we have
	
	\begin{equation}{\label{largeq}}
		g^{(k)}(q,x) \ll_{k,\epsilon} Q^{ \epsilon k}
	\end{equation}
	
	which says that for $q\gg Q^{1-\epsilon}$, $g(q,x)$ has essentially no oscillation. To prove this, by definition we have
	
	$$g(q,x) = \int \Delta_q(u)f(u)e \left( - \frac{ux}{qQ}\right)du$$
	
	where $\Delta_q(u)$ satisfies
	
	\begin{equation}\label{Deltabound}
		\Delta_q(u) \ll \frac{1}{Q(q+Q)} + \frac{1}{|u| + qQ}
	\end{equation}
	
	and $f$ is a smooth function supported in $[-Q^2,Q^2]$ with $f(0) = ||f||_{L^\infty} = 1$ and $f^j(u) \ll Q^{-2j}$ (see \cite{iwaniec_kowalski})
	
	Then by differentiating under the integral sign we have
	
	\begin{equation}\label{defng(q,x)}
		g^{(k)}(q,x) = \left(-\frac{2 \pi i}{qQ}\right)^k \int \Delta_q(u)f(u)u^k e \left( - \frac{ux}{qQ}\right)du
	\end{equation}
	
	Applying (\ref{Deltabound}), (\ref{defng(q,x)}) and assuming $q \gg Q^{1-\epsilon}$ we have
	
	\begin{equation}
		\begin{aligned}
			g^{(k)}(q,x)  & \ll \frac{Q^{\epsilon k}}{Q^{2(k+1)}}\left( \int \limits_{|u| \ll Q^{2}} |u|^k du + Q^2\int \limits_{|u| \ll Q^2} |u|^{k-1} du \right )
			\\ & \ll Q^{\epsilon k}
		\end{aligned}
	\end{equation}
	
	Applying the delta symbol (\ref{smallx}) we have
	
	\begin{equation}{\label{DFI2}}
		\begin{aligned}
			\delta(m=n^2) =&  \sum \limits_{q\leq Q} \frac{1}{Q} \sideset{}{^*}\sum \limits_{a \bmod q} \frac{1}{q} e \left( \frac{a(n^2-m)}{q} \right) \int \limits_\R W(x)g(q,x) e \left( \frac{(n^2-m)x}{qQ} \right) dx
			\\ &+ O (Q^{-A})
		\end{aligned}
	\end{equation}
	
	with $W(x)$ a smooth bump function supported in $[-2Q^\epsilon,2Q^\epsilon]$ with $W(x) \equiv 1$ on $[-Q^\epsilon,Q^\epsilon]$ and satisfying $W^{(i)}(x) \ll_i 1$.
	
	Applying $(\ref{DFI2})$ to (\ref{defnS(N)}) we get
	
	\begin{equation}
		\begin{aligned}
			S(N) =&  \sum \limits_{q \leq Q} \int \limits_{\R} W(x) \frac{g(q,x)}{qQ} \sideset{}{^*}\sum \limits_{a \bmod q} \sum \limits_m \lambda_f(m) U \left( \frac{m}{N^2} \right) e\left( -\frac{am}{q}- \frac{mx}{qQ}\right)
			\\ & \times \sum \limits_n V \left( \frac{n}{N}\right) n^{-it}e \left( \frac{an^2}{q} + \frac{n^2x}{qQ}\right)dx dv
		\end{aligned}
	\end{equation}
	
	plus a negligible error term. Here $U(x)$ is another bump function supported on $[3/4, 9/4]$ and $U(x) \equiv 1$ on $[1,2]$ satisfying the bounds $U^{(i)}(x) \ll_i 1$.
	
	\begin{thm}\label{thm2}
		
		For $pt \ll N \ll p^{1+\epsilon}t^{3/2+\epsilon}$ and $t^{\epsilon} \ll K \ll t^{1-\epsilon}$ we have
		
		$$S(N) \ll N^{1/2}p^{1/2+\epsilon} t^{1/2+\epsilon} \left( K^{1/4} + \frac{t}{K^{5/4}}\right)$$
		
	\end{thm}
	
	Choosing the optimal choice of $K = t^{2/3}$ gives Theorem \ref{thm1}.

	\subsection{Outline}
	
	For simplicity let's assume that we're in the generic situation i.e. $N = pt^{3/2}$ and $q \sim Q$. dropping the $\epsilon$s. Then we're looking at an expression (roughly) of the form
	
	$$ \sum \limits_{q \sim Q}~\sideset{}{^*}\sum \limits_{a \bmod q} \sum \limits_{m \sim N^2} \lambda_f(m) e \left(-\frac{am}{q}\right) \sum \limits_{n \sim N} n^{-it} e \left( \frac{an^2}{q}\right)$$
	
	where we're ignoring the analytic oscilation coming from the integral over $x$. Notice that by applying the delta method we've got $p^3t^{9/2}$ terms. So to beat the trivial convexity estimate we need to save $p^2t^3$ and a little more.
	\\
	
	We apply Voronoi summation to the $m$ sum which gives us a sum of length $m \sim pK$ which saves us $N/\sqrt{pK}$. We apply Poisson summation to the $n$ sum which gives us a sum of length $n \sim t/\sqrt{K}$ which saves us $\sqrt{N/(t/\sqrt{K})}$. The integral over $x$ saves us an extra $\sqrt{K}$ and the resulting character sum arising from the sum over $a \bmod q$ saves us another $\sqrt{q} \sim \sqrt{Q}$. So far we have saved
	
	$$\frac{N}{\sqrt{pK}} \times \sqrt{ \frac{N}{t/\sqrt{K}}} \times \sqrt{K} \times \sqrt{Q} = \frac{N}{\sqrt{pt}}$$
	
	It remains to save $\sqrt{pt}$.
	
	We've arrived at the expression (after simplifying the sum $a \bmod q$)
	
	$$\sum \limits_{q \leq Q} \int \limits_x \frac{g(q,x)}{qQ} \sum \limits_{m \sim pK} \sum \limits_{n \sim t/\sqrt{K}} \lambda_f(m) \left( \frac{4m - pn^2}{q}\right) \int \limits_v (...)$$
	
	Applying Cauchy-Schwarz in the $q$ variable,we get the following expression
	
	$$L = \sum \limits_{q \sim Q} \left| \sum \limits_{d \sim pt^2/\sqrt{K}} \left( \frac{d}{q} \right) \alpha(d)\right|^2$$
	
	where we've set $d = pn^2 - 4m$ and
	
	$$\alpha(d) = \underset{\substack{m \sim pK\\ n \sim t/\sqrt{K} \\ d = pn^2-4m}} {\sum \sum} \lambda_f(m)\int \limits_v (...) \ll (pt)^\epsilon \underset{\substack{m \sim pK\\ n \sim t/\sqrt{K} \\ d = pn^2-4m}} {\sum \sum} 1$$
	
	Since we applied Cauchy-Schwarz we now require a saving of $pt$. Trivially, we see that
	
	$$L \ll Q(pK)^2(t/\sqrt{K})^2 =  p^3 t^{7/2}\sqrt{K}$$
	
	We now apply quadratic large sieve which gives (ignoring the square-free conditions on $q$ and $d$)
	
	$$L \ll (pt)^\epsilon \left(Q + pt^2/K\right) ||\alpha||^2$$
	
	where
	
	$$||\alpha||^2 = \sum \limits_d |\alpha(d)|^2 = \underset{\substack{m_1,m_2 \sim pK\\ n_1,n_2 \sim t/\sqrt{K} \\ pn_1^2 - 4m_1 = pn_2^2-4m_2}} {\sum \sum \sum \sum} 1$$
	
	We get $n_1^2 - n_2^2 =4(m_2-m_1)/p \ll K$ which implies $n_1-n_2 \ll K/(n_1+n_2) \ll K^{3/2}/t$. If $t^{2/3} \ll K \ll t$ then the number of pairs of $(n_1,n_2)$ is $\ll K$. It follows that the total number of terms  $(n_1,n_2,m_1,m_2)$ is $\ll pK^2$. So $||\alpha||^2 \ll pK^2$ and hence $L \ll (pt)^\epsilon p^2t^2K$. This saves us $pt^{3/2}/\sqrt{K}$.  This is better than $pt$ as long as $K < t$. Choosing $K = t^{2/3}$ gives the best possible result.
	
	\subsection{Quadratic Gauss Sums}
	
	Here we state the evaluation of certain quadratic Gauss sums. Define
	
	$$G(a,b;c) = \sum \limits_{x \bmod c} e \left( \frac{ax^2}{c} + \frac{bx}{c}\right)$$
	
	We always assume $(a,c)=1$.
	
	\begin{lem}\label{qgauss1}
		
		(1) If $c$ is even, say $c=2r$ then $G(a,b,2r)$ is $0$ unless $ar+b$ is even.
		
		(2) If $(r,s)=1$ then $G(a,b,rs)=G(ar,b,s)G(as,b,r)$.
		
		(3) We have
		
		\begin{equation}G(a,b,c) = \begin{cases}
				e\left(-\frac{\overline{4a}b^2}{c}\right)G(a,0,c) & \text{if } c \text{ is odd}
				\\ 2e \left(- \frac{\overline{8a_r}b^2}{r}\right)G(2a,0,r) & \text{if } 2||c,c=2r \text{ and } b \text{ is odd}
				\\
				e \left(-\frac{ \overline{a}(b/2)^2}{c}\right)G(a,0,c) & \text{if } 4|c \text{ and } b \text{ is even}
			\end{cases}
		\end{equation}
		
		$a_r$ is to denote the inverse mod $r$.
	\end{lem}
	
	\textit{Proof}: This is contained in Lemma 5.4.5 of \cite{huxley}.
	\qed
	\begin{lem}\label{qgauss2}
		
		Assume $c$ is odd and let $(a,c)=1$. We have
		
		\begin{equation}\frac{G(a,0,c)}{c^{1/2}} = \begin{cases}
				\left( \frac{a}{c}\right) & \text{if } c \equiv 1 \bmod 4
				\\i\left( \frac{a}{c}\right) & \text{if } c \equiv 3 \bmod 4
			\end{cases}
		\end{equation}
	\end{lem}
	
	\textit{Proof}: This is Theorem 1.3.4 of \cite{gauss_jacobi_sums}.
	\qed
	\\
	We denote
	
	$$\epsilon_c = \begin{cases} 1 & \text{if } c \equiv 1 \bmod 4 \\ i & \text{if } c \equiv 3 \bmod 4 \end{cases}$$
	
	\begin{lem}\label{qgauss3}
		
		Suppose $c = 2^kc_0$ with $k\geq 2$ and $c_0$ odd. Assume $(a,c)=1$. Then
		
		\begin{equation}\frac{G(a,0,c)}{c^{1/2}} = \epsilon_{c_0}\left( \frac{a}{c_0}\right)\begin{cases} 1+e\left( \frac{ac_0}{4}\right) & \text{if } k \equiv 0 \bmod 2
				\\ 2^{1/2}e\left( \frac{ac_0}{8}\right) & \text{if } k \equiv 1 \bmod 2
			\end{cases}
		\end{equation}
		
	\end{lem}
	
	\textit{Proof}: This is Lemma 5.1 of \cite{khan_moments_2021}.
	\qed
	\\
	\subsection{Voronoi Summation}
	
	Now we state the Voronoi summation formula.
	
	\begin{thm}[Kowalski-Michel-VanderKam]
		Let $g(x)$ be a smooth function on $\R_{> 0}$ with compact support. Let $(a,q) = 1$ and let $\overline{a}$ denote the multiplicative inverse of $a$ mod $q$. Let $p$ be a prime and assume $f$ is a newform of level $p$. Set $p_2 = p/(p,q)$. Then there exists a newform $f^*$ of the same level and weight and a complex number $\eta$ of modulus 1 (which depends on $f,p$) s.t.
		
		$$\sum \limits_{n=1}^\infty e \left( \frac{an}{q}\right) \lambda_f(n)g(n) = \frac{2 \pi \eta}{q\sqrt{p_2}} \sum \limits_n \lambda_{f^*}(n) e \left( - \frac{n \overline{ap_2}}{q}\right) \int \limits_0^\infty h(y) J_{k-1} \left( \frac{4 \pi \sqrt{ny}}{q \sqrt{p_2}} \right) dy$$
	\end{thm}
	
	\textit{Proof}: See the appendix to  \cite{kowalski_rankin-selberg_2002}.
	\qed
	\\
	\subsection{Stationary Phase Analysis}
	
	In the proof, we require a study of oscillatory integrals of the form
	
	$$\int g(x)e(f(x))dx$$
	
	with $g$ and $f$ functions on $\R$. We use the following less precise version of stationary phase analysis due to Mckee, Sun and Ye.
	
	\begin{thm}[Mckee-Sun-Ye]\label{thmsphase}
		
		Let $f(x)$ and $g(x)$ be real-valued smooth functions on $\R$. Define
		
		$$H_1(x)  = \frac{g(x)}{2 \pi i f^\prime(x)}$$
		$$H_i(x) = \frac{H_{i-1}^\prime(x)}{2\pi i f^{\prime(x)}}$$
		for $i \geq 2$.
		
		Let $[\alpha,\beta]$ be an interval. Assume there are parameters $M,N,T,U$ s.t.
		
		$$M \geq \beta - \alpha$$
		
		and positive constants $C_k$ s.t. for $\alpha \leq x \leq \beta$ s.t.
		
		$$|f^{(r)}(x)| \leq C_r\frac{T}{M^r}~~~r \geq 2$$
		$$f^{\prime \prime}(x) \geq \frac{T}{C_2M^2}$$
		
		and
		
		$$|g^{(s)}(x)| \leq C_s\frac{U}{N^s}~~~s \geq 0$$
		
		Suppose $f^\prime(x)$ changes sign only at $x = \gamma$ from negative to positive with $\alpha < \gamma < \beta$. Let
		
		$$\Delta_n = \inf \left\{ \frac{\log 2}{C_2}, \frac{1}{C_2^2 \sup \limits_{2\leq k \leq 2n+3}\{C_k\}}\right\}$$
		
		If $T\Delta_n > 1$ then for $n \geq 2$ we have
		
		\begin{equation}
			\begin{aligned}
				&\int \limits_\alpha^\beta g(x)e(f(x))dx
				\\ = & \frac{e(f(\gamma)+1/8)}{\sqrt{f^{\prime \prime}(\gamma)}}\left( g(\gamma) + \sum \limits_{j=1}^n \omega_{2j} \frac{(-1)^j(2j-1)!}{(4\pi i\lambda_2)^j}\right) + \left[ e(f(x))\sum \limits_{i=1}^{n+1}H_i(x)\right]_\alpha^\beta
				\\ & + O\left(\frac{UM^{2n+5}}{T^{n+2}N^{n+2}} \left( \frac{1}{(\gamma - \alpha)^{n+2}} + \frac{1}{(\beta - \gamma)^{n+2}}\right)\right)
				\\ & + O \left( \frac{UM^{2n+4}}{T^{n+2}}\left( \frac{1}{(\gamma-\alpha)^{2n+3}} + \frac{1}{(\beta-\gamma)^{2n+3}}\right)\right)
				\\ & + O \left(\frac{UM^{2n+4}}{T^{n+2}N^{2n}} \left( \frac{1}{(\gamma-\alpha)^3} + \frac{1}{(\beta-\gamma)^3}\right)\right) + O \left( \frac{U}{T^{n+1}}\left( \frac{M^{2n+2}}{N^{2n+1}} + M\right)\right)
			\end{aligned}
		\end{equation}
		
		where
		
		$$\lambda_k = \frac{f^{(k)}(\gamma)}{k!}~~~k \geq 2$$
		$$\eta_k = \frac{g^{(k)}(\gamma)}{k!}~~~k \geq 1$$
		
		and
		
		$$\omega_k = \eta_k + \sum \limits_{l=0}^{k-1}\eta_l\sum \limits_{j=1}^{k-l} \frac{C_{k,l,j}}{\lambda_2^j} \sum \limits_{\substack{3 \leq n_1,...,n_j \leq 2n+3 \\ n_1+n_2+...+n_j = k-l+2j}} \lambda_{n_1}\lambda_{n_2}...\lambda_{n_j}$$
		
		with $C_{k,l,j}$ constants.
	\end{thm}
	
	\textit{Proof}: This is Theorem 1.2 of \cite{mckee_weighted_2017}.
	\qed
	\\
	
	We've stated the theorem only when $f(x) = x^{-it}e(rx)$ with $r \in \R$ and $g(x)$ is a smooth function with compact support. The original theorem is stated in more generality. We will choose $\alpha$ and $\beta$ away from the support of $g(x)$ so that $H_i(\alpha)=H_i(\beta) = 0$ for all $i$.
	
	Had we only been interested in just the $t$-aspect, we would only need to expand upto $n = 2$ with a sufficient error term. However, due to the presence of the level in the $GL(2)$ form, we require more terms in the expansion so that the error term is much smaller.
	
	\subsection{Quadratic Large Sieve}
	
	We state a version of the large sieve when the characters run over only quadratic characters, a result due to Heath-Brown.
	
	\begin{thm}[Heath-Brown]{\label{qsieve}}
		
		Let $Q,N \geq 1$ and let $a_n$ be a sequence of complex numbers for $1 \leq n \leq N$. Then
		
		$$\sideset{}{^*}\sum \limits_{q \leq Q} \left|\sum \limits_{n \leq N} a_n \left(\frac{n}{q} \right) \right|^2 \ll (QN)^\epsilon(Q+N)\sum \limits_{n_1n_2 = \square} |a_{n_1}a_{n_2}|$$
		
		where $\sideset{}{^*}\sum$ indicates that the sum is taken over square-free integers.
	\end{thm}
	
	\textit{Proof}: This is Corollary 2 of \cite{quadratic_large_sieve}.
	\qed
	\subsection{Mellin Transforms}
	
	Suppose $V(x)$ is a smooth function with compact support in $\R_{>0}$. Then we have that
	
	$$V(x) = \int \limits_{c-i\infty}^{c+i\infty} \tilde{V}(s) x^{-s} ds = \int \limits_{(c)} \tilde{V}(s) x^{-s} ds$$
	
	where $c > 0$ and $\tilde{V}(s)$ is the Melling transform of $V$, given by
	
	$$\tilde{V}(s) = \int \limits_0^\infty V(x)x^{s-1}dx$$
	
	If $V(x)$ has compact support(or even Schwarz), then it's Mellin transform has rapid decay in vertical strips i.e. $s = \sigma + i\tau$ with $\sigma \in [\sigma_0,\sigma_1]$ with $0 < \sigma_0 < \sigma_1$ we have
	
	$$\tilde{V}(s) \ll_{A,\sigma_0,\sigma_1} (1 + \tau)^{-A}$$
	
	for all $A\geq 1$.
	
	\section{Application of Summation Formulae}
	
	\subsection{Poisson Summation}
	
	Here we apply the Poisson summation to the $n$ sum. We set $g(y) = y^{it - 2iv} V \left( \frac{y}{N} \right)e \left( \frac{y^2 x}{qQ}\right)$.
	
	\begin{equation}\begin{aligned}
			\sum \limits_{n\geq 1} g(n) e \left( \frac{an^2}{q} \right) & = \sum \limits_{r \bmod q} e \left(\frac{ar^2}{q} \right) \sum \limits_{n \equiv r \bmod q} g(n)
			\\ & = \frac{1}{q}\sum \limits_{n\in \Z} G(a,n,q)\int \limits_\R y_1^{-it- 2iv} V \left( \frac{y_1}{N}\right)e \left( \frac{y_1^2x}{qQ} - \frac{ny_1}{q}\right) dy_1
		\end{aligned}
	\end{equation}
	
	Making the change of variable $y \mapsto yN$ we see that by repeated integration by parts that $I(n,q)$ is negligibly small if
	
	\begin{equation}n \gg \left( \frac{qt}{N} + K^{1/2}\right) = :N_0
	\end{equation}
	
	After applying the change of variables we get
	
	$$\frac{N^{1-it}}{q}\sum \limits_{|n| \ll N_0} G(a,n,q) I_1(n,q)$$
	
	where
	
	$$I_1(n,q) = \int \limits_\R y_1^{-it} V \left(y_1\right)e \left( \frac{N^2y_1^2 x}{qQ} - \frac{nNy_1}{q}\right) dy_1$$
	
	\subsection{Voronoi Summation}
	
	Here we set our $g(y) = y^{iv}U( \frac{y}{N^2})e( -\frac{xy}{qQ})$. Then by Voronoi summation
	
	\begin{equation}\begin{aligned}
			&\sum \limits_{m=1}^\infty \lambda_f(m)m^{iv} U \left( \frac{m}{N^2} \right) e\left( -\frac{am}{q}- \frac{mx}{qQ}\right)
			\\ =& \frac{1}{q\sqrt{p_2}} \sum \limits_{m=1}^\infty \lambda_{f^*}(m) e \left( \frac{m\overline{ap_2}}{q}\right) \int \limits_\R U \left( \frac{y_2}{N^2} \right) e \left( \frac{y_2x}{qQ} \right) J_{k-1}\left( \frac{4 \pi \sqrt{my_2}}{q\sqrt{p_2}}\right)dy_2
	\end{aligned}\end{equation}
	
	We have the identity $2J_{k-1}(x) = e^{ix} W_{k-1}^+(x) + e^{-ix}W_{k-1}^-(x)$ with
	
	\begin{equation}{\label{jbound}}
		x^j W_{k-1}^{(j)}(x) \ll_{k,j} (1+x)^{-1/2}\end{equation}
	
	Therefore we may write the integral as a linear combination of integrals of the form
	
	\begin{equation}\int \limits_\R  U \left(\frac{y_2}{N^2}\right) e \left(- \frac{y_2 x}{qQ} \right) e \left(\pm \frac{2\sqrt{my_2}}{q\sqrt{p_2}}\right) dy_2
	\end{equation}
	
	Recall that $p_2 = p/(p,q)$. We from now on assume $p_2 = p$ i.e. $(p,q)=1$. The case $p|q$ gives a considerably shorter sum and saves an extra $\sqrt{p}$ over the complimentary case.
	
	Let us make the change of variable $y \mapsto yN^2$ to obtain
	
	$$I_2(m,q) = \sum \limits_{\pm} \int \limits_\R U(y_2)e \left( - \frac{N^2y_2x}{qQ} \pm \frac{2N\sqrt{my_2}}{q\sqrt{p}}\right) dy_2$$
	
	and by repeated integration by parts we have that the integral is negligible if
	
	\begin{equation}\label{msize}
		m \gg N^\epsilon pK =:M_0
	\end{equation}
	
	We get that the $m$-sum is transformed into
	
	$$\frac{N^{2+2iv}}{p^{1/2}q} \sum \limits_{m \ll M_0} \lambda_f(m) W_{k-1} \left( \frac{N\sqrt{m}}{q\sqrt{p}}\right) e \left( \frac{\overline{ap}m}{q}\right) I_2(m,q)$$

	\section{Simplification}

	After applying the summation formulae we arrive at (a linear combination of)
	
	$$\frac{N^{3-it}}{QKp^{1/2}}\sum \limits_{\substack{q \leq Q \\ (q,p)=1}} \frac{1}{q^2}\sum \limits_{|n| \ll N_0}\sum \limits_{m \ll M_0} \lambda_f(m)W_{k-1} \left(\frac{2N\sqrt{m}}{q\sqrt{p}}\right) \mathcal{C}(m,n,q)I(m,n,q)$$
	
	where $\mathcal{C}(m,n,q)$ is a character sum
	
	\begin{equation}\mathcal{C}(m,n,q) = \frac{1}{q}\sideset{}{^*}\sum \limits_{a \bmod q} G(a,n,q)e\left( \frac{\overline{ap}m}{q}\right)
	\end{equation}
	
	and $I(m,n,q)$ is a threefold integral
	
	\begin{equation}{\label{I}}\begin{aligned}
			I(m,n,q) = & \int \limits_\R W(x)g(q,x)\int \limits_\R U(y_1)y_1^{-it} \int \limits_0^\infty V(y_2)
			\\ & \times e \left( \frac{N^2x(y_1^2-y_2)}{qQ} - \frac{nNy_1}{q} -\pm \frac{2N\sqrt{my_2}}{qp^{1/2}}\right) dy_1dy_2dx
		\end{aligned}
	\end{equation}
	
	We now need to explicitly calculate the character sum and simplify the integral.
	
	\subsection{The Character Sum}
	
	First we evaluate the character sum $\mathcal{C}(m,n,q)$.
	
	We first assume that $q$ is odd. Applying Lemma (\ref{qgauss1}) and Lemma (\ref{qgauss2}) and a change of variable $a \mapsto \overline{a} \bmod q$ we get
	
	$$\frac{1}{q^{1/2}}\sideset{}{^*}\sum \limits_{a \bmod q} \left( \frac{a}{q}\right)e \left( \frac{(\overline{p}m - \overline{4}n^2)a}{q}\right)$$
	
	Next make the change of variable $a \mapsto 4ap$ giving
	
	$$\frac{1}{q^{1/2}}\left(\frac{p}{q}\right) \sum \limits_{a \bmod q} \left( \frac{a}{q}\right)e \left( \frac{a(4m-pn^2)}{q}\right)$$
	
	Denote the sum over $a \bmod q$ as $\mathcal{C}_1(4m-pn^2,q)$. For notational purposes, let's assume $l=4m - pn^2$.

	First, decompose $q=q^*q_0^2$ with $q^*$ square-free. Then $q^*$ is the conductor of the Legendre symbol $\left( \frac{\cdot}{q}\right)$. Next, we split $q_0^2 = q_1q_2$ where $q_1|(q^*)^\infty$ and $(q_2,q^*)=1$. We apply Chinese Remainer Theorem getting
	
	$$\begin{aligned}
		\mathcal{C}_1(l,q) & = \sum \limits_{a \bmod q_1q^*} \left( \frac{a}{q^*} \right) e \left( \frac{a \overline{q_2}l}{q_0}\right) \sideset{}{^*}\sum \limits_{b \bmod q_2} e \left( \frac{b \overline{q_0}l}{q_2}\right)
		\\ & = S(l,0,q_2) \sum \limits_{a \bmod q_1q^*} \left( \frac{a}{q^*}\right)e \left( \frac{al}{q_0}\right)
		\\ & = S(l,0,q_2) \sum \limits_{a \bmod q^*} \left( \frac{a}{q^*}\right) \sum \limits_{k=1}^{q_1} e \left( \frac{l(a+kq^*)}{q_0}\right)
		\\ & = \begin{cases}
			S(l,0,q_2)q_1\epsilon_{q^*}\sqrt{q^*}\left( \frac{\frac{l}{q_1}}{q^*}\right) & \text{if } q_1|l
			\\ 0 & \text{otherwise}
		\end{cases}
	\end{aligned}$$
	
	Next we assume $2||q$, say $q=2q_0$ with $q_0$ odd. Then the character sum is given by
	
	$$\frac{1}{q_0^{1/2}}\sideset{}{^*}\sum \limits_{a \bmod 2q_0} \left( \frac{2a}{r}\right)e \left( \frac{\overline{ap}m}{2q_0} - \frac{\overline{8a_{q_0}}n^2}{q_0}\right)$$
	
	where $a_{q_0}$ denote the residue class of $a \bmod q_0$. Then we note that $\overline{a} \equiv 2 \overline{2a_r} +r \bmod 2r$. Applying a suitable change of variables gives
	
	$$\frac{2}{q_0^{1/2}}\left(\frac{p}{q_0}\right)e\left( \frac{pm}{2}\right)\mathcal{C}_1(4m - pn^2,q_0)$$
	
	The character sum $\mathcal{C}_1$ can then be evaluated as above.
	
	Now assume $4|q$. Suppose $q=2^kq_0$ with $q_0$ odd and $k \geq 2$. First we assume $k \equiv 0 \bmod 2$. Then by Lemma \ref{qgauss3} we have
	
	$$q^{1/2}\epsilon_{q_0} \sideset{}{^*}\sum \limits_{a \bmod q} \left( \frac{a}{q_0}\right) e \left( \frac{(\overline{p}m - (n/2)^2)\overline{a}}{q}\right)\left(1 + e \left( \frac{aq_0}{4}\right) \right)$$
	
	We apply a change of variable $a \mapsto ap$ to get
	
	$$q^{1/2}\epsilon_{q_0} \left( \frac{p}{q_0}\right) \sideset{}{^*}\sum \limits_{a \bmod q} \left( \frac{a}{q_0}\right) e \left( \frac{(m - p(n/2)^2)\overline{a}}{q}\right)\left(1 + e \left( \frac{apq_0}{4}\right) \right)$$
	
	Next we make the change of variable $a \mapsto \overline{a} \bmod c$ and note that $a \equiv \overline{a}_c \bmod 4$. Thus we get
	
	$$q^{1/2}\epsilon_{q_0}\left( \frac{p}{q_0}\right)\sideset{}{^*}\sum \limits_{a \bmod q} \left( \frac{a}{q_0}\right)
	e \left( \frac{(m - p(n/2)^2)a}{q}\right)\left(1 + e \left( \frac{apq_0}{4}\right)\right)$$
	
	Applying CRT gives
	
	$$\begin{aligned}
		q^{1/2}\epsilon_{q_0} \left( \frac{p}{q_0}\right) \sum \limits_{a \bmod q_0} \left( \frac{a}{q_0}\right)e \left( \frac{\overline{2^k}a(m - p (n/2)^2)}{q_0}\right)
		\\ \times \sideset{}{^*}\sum \limits_{b \bmod 2^k} e \left( \frac{b\overline{q_0}(m - p (n/2)^2)}{2^k}\right)\left(1 + e \left( \frac{bpq_0}{4}\right)\right)
	\end{aligned}$$
	
	We now calculate the $b \bmod 2^k$ sum. For notation, we set $l = m -p(n/2)^2$. We have
	
	$$\begin{aligned}
		\sideset{}{^*}\sum \limits_{b \bmod 2^k} e \left( \frac{b \overline{q_0}l}{2^k}\right)\left( 1 + e \left( \frac{bpq_0}{4}\right) \right) & = \sideset{}{^*}\sum \limits_{b \bmod 2^k} e \left( \frac{bl}{2^k}\right)\left( 1 + e \left( \frac{bp}{4}\right)\right)
		\\ & = S(l,0,2^k) + \sideset{}{^*}\sum \limits_{b \bmod 2^k} e \left( \frac{b l}{2^k}\right)e \left( \frac{bp}{4}\right)
		\\ & = S(l,0,2^k) + \sideset{}{^*}\sum \limits_{a \bmod 4}e \left( \frac{ap}{4}\right)\sum \limits_{r=1}^{2^{k-2}}e \left( \frac{(a + 4r)l}{2^k}\right)
		\\ & = S(l,0,2^k) + S(p+l,0,4)2^{k-2}\delta(2^{k-2}|l)
	\end{aligned}$$
	
	Note that $S(p+l,0,4)=0$ if $2|l$. Hence the second term vanishes unless $k=2$.
	
	Similarly evaluating in the case $k \equiv 1 \bmod 2$ and noting that $a \equiv \overline{a_q} \bmod 8$, we get
	
	$$\mathcal{C}_1(l,q_0)S(p+l,0,8)2^{k-3}\delta(2^{k-3}|l)$$
	
	Similarly, $S(p+l,0,8) = 0$ when $2|l$. Hence there is contribution only when $k=3$.
	
	We assume that $q$ is odd, since the even case has an additional Ramnujan sum or a constant weight factor, either of which provides little additional complication.
	
	\subsection{The Integral}
	
	Now we attempt to simplify the threefold integral. We follow the treatment in \cite{munshi_subconvexity_2021}. Recall this is given by $(\ref{I})$
	
	$$\begin{aligned}
		I(m,n,q) = & \int \limits_\R W(x)g(q,x)\int \limits_\R U(y_1)y_1^{-it} \int \limits_0^\infty V(y_2)
		\\ & \times e \left( \frac{N^2x(y_1^2-y_2)}{qQ} - \frac{nNy_1}{q} \pm \frac{2N\sqrt{my_2}}{qp^{1/2}}\right) dy_1dy_2
	\end{aligned}$$
	
	Now we assume that $q \sim C$ and $m \sim M_1$ with $C \ll Q$ and $M_1 \ll M$. First we assume that $N^2M/C^2p \gg N^\epsilon$. Then we note that looking at the $y_2$ integral $I_2(m,q)$ we see that $I(m,n,q)$ is negligible unless $t^{-\epsilon}\sqrt{\frac{M_1}{pK}} \ll |x| \ll t^\epsilon \sqrt{ \frac{M_1}{pK}}$.
	
	Next assume that $M_1 \ll C^2p/N^{2-\epsilon}$. Then looking at the $y_2$ integral once again, we see that the integral is negligible unless $|x| \ll CQ/N^{2-\epsilon}$.
	
	Therefore we may control the domain of integration of the $x$-integral as a function of $M_1$, say $I(M_1)$ with
	
	\begin{equation}
		\mu(I(M_1)) \ll \begin{cases} CQ/N^{2-\epsilon} & \text{if }M_1 \ll C^2p/N^{2-\epsilon}
			\\ t^\epsilon\sqrt{\frac{M_1}{pK}} & \text{if } M_1 \gg C^2p/N^{2-\epsilon}
		\end{cases}
	\end{equation}
	
	Now we assume $C \ll Q^{1-\epsilon}$. By looking at the $x$-integral and using (\ref{smallq}) we conclude that $|y_1^2-y_2| \ll N^{2\epsilon}C/QK$ which implies $|y_1 - \sqrt{y_2}|\ll N^{2\epsilon}C/QK$. Next, assume $q \gg Q^{1-\epsilon}$. Then we use (\ref{largeq}) and repeated integration by parts to conclude that $|y_1 - \sqrt{y_2}|\ll N^{2\epsilon}C/QK$. So for all $q \leq Q$ we conclude that $|y_1 - \sqrt{y_2}|\ll N^{2\epsilon}C/QK$.
	
	We set $\sqrt{y_2}=y_1+u$ with $|u|\ll N^{2\epsilon}C/QK$ and insert this into the integral to obtain
	
	$$\begin{aligned}
		&\int \limits_{\R} W(x)g(q,x)\int \limits_\R U(y_1)y_1^{-it} \int \limits_0^\infty V((y_1 + u)^2)\sqrt{y_1+u}
		\\ & \times e \left( \frac{N^2x(y_1^2-(y_1+u)^2)}{qQ} - \frac{nNy_1}{q} \pm \frac{2N\sqrt{m}(y_1+u)}{q\sqrt{p}}\right)dudy_1dx
	\end{aligned}$$

	Now we want to perform stationary phase analysis on the $y_1$ integral. The weight function is given by
	
	$$G(z) = U(z)V((z+u)^2)\sqrt{z+u}$$
	
	The phase function is given by
	
	$$P(z) = -\frac{t}{2\pi}\log z - \frac{2N^2uxz}{qQ} - \frac{nNz}{q} \pm \frac{2N\sqrt{m}z}{q\sqrt{p}}$$
	
	$$P^\prime(z) = -\frac{t}{2\pi z} - \frac{2N^2ux}{qQ} - \frac{nN}{q} \pm \frac{2N\sqrt{m}}{q\sqrt{p}}$$
	
	$$P^{(j)}(z) = t\frac{(-1)^j(j-1)!}{2\pi z^j}~~~ j \geq 2$$
	
	Therefore we see that the stationary point is given by
	
	$$z_0 = \frac{2\pi N}{qt}p(m,n,u,x)^{-1}$$
	
	where
	
	$$p(m,n,u,x) = \frac{2Nux}{Q} -n \pm \sqrt{\frac{m}{p}}$$
	
	For notational purposes, we'll drop the dependence on $u$ and $x$. We apply Theorem \ref{thmsphase} with $\alpha = 1/8$, $\beta = 3$, $M = \beta- \alpha$, $T = t$ and $U=N = 1$.
	
	Due to support conditions on the weight function, we must have $z_0 \asymp 1$ otherwise the integral is negligibly small. This implies
	
	\begin{equation}p(m,n) \asymp \frac{Ct}{N}\end{equation}
	
	Note that
	
	\begin{equation}\label{smallshift}
		\frac{2Nux}{Q} \ll \frac{1}{N}\end{equation}
	
	If in addition, $qt/N \gg t^\epsilon K^{1/2}$ then $n \sim Ct/N$ from (\ref{smallshift}) and  (\ref{msize}).
	
	\subsubsection{Separation of Variables}
	
	Before, we apply the quadratic large sieve, we need to make sure that we have a bilinear form in the $q$ and $m,n$ variables. Towards this, let's make some observations.
	
	First, $P^{\prime \prime}(z_0) \sim t$. Second, $P(z_0)$ is a function independent of $q$ (but will still depend on $C$). Third, $P^{(j)}(z_0)$ can be written as $a(q)b(m,n)$ with $a(q)$ and $b(m,n)$ are two functions. So we see that the only place where the integral is not separated is in the weight function $G(z_0)$ and $W_{k-1}\left( \frac{4\pi N\sqrt{m}}{q\sqrt{p}}\right)$.
	
	To handle this, we apply Mellin inversion. Let
	
	$$\tilde{G}(s) = \int \limits_0^\infty G(x)x^{s-1}dx$$
	
	be the Mellin transform of $G$ valid for $\Re(s) > 0$. Then by Mellin inversion
	
	$$G(z_0) = \frac{1}{2\pi i }\int \limits_{(c)} \tilde{G}(s) z_0^s ds$$
	
	valid for $c > 0$. Since $\tilde{G}(s)$ has rapid decay in vertical strips, we may truncate the integral upto height $(pt)^\epsilon$ at the cost of a negligible error term. We also set $c=\epsilon$. This allows us to separate the parameters at a cost of a $O((pt)^\epsilon)$ multiplier. The same is true for all of it's derivatives (constants depending on the order of the derivative taken).
	
	Similarly, we apply Mellin inversion to get
	
	$$W_{k-1} \left( \frac{4 \pi N\sqrt{m}}{q\sqrt{p}}\right) = \frac{1}{2\pi i} \int \limits_{(c)} \tilde{W}_{k-1}(s) \left( \frac{4\pi N \sqrt{m}}{q\sqrt{p}} \right)^{s}ds$$
	
	where
	
	$$\tilde{W}_{k-1}(s) = \int \limits_0^\infty W_{k-1}(x)x^{s-1}dx$$
	
	By (\ref{jbound}), we may replace $W_{k-1} \left( \frac{4\pi N\sqrt{m}}{q\sqrt{p}} \right)$ by $N^\epsilon \left(1 + \frac{4\pi N \sqrt{m}}{q\sqrt{p}} \right)^{-1/2}$.
	Now we make some notational simplifications. First, we assume that we've performed all the necessary variable separations, and drop the integrals over $s$. Second, we assume that $P^{\prime \prime}(z_0) = ct$ with $c > 0$ a fixed constant. This simplifies further calculations without altering the final boudn. Third, we only make calculations with $G(z)$ as the lower order terms in the stationary phase expansion can be dealt with similarly but with better bounds.
	
	Finally, we observe that $u$ depends on the size of $q$, which we may handle by applying a dyadic partition on the $q$-variable. This gives
	
	\begin{equation}{\label{qdyadic}}
		S(N) = \frac{N^{3-it}}{Qp^{1/2}}\sum \limits_{C \text{ dyadic}} S(N;C) + O(N^{-A})
	\end{equation}
	
	where
	
	$$\begin{aligned}
		S(N,C) =& \sum \limits_{q\sim C} \frac{1}{q^{3+it}}
		\\ & \times \sum \limits_{n \ll N_0}\sum \limits_{m \ll M_0}\lambda_f(m) W_{k-1} \left( \frac{4 \pi N\sqrt{m}}{q\sqrt{p}} \right)\mathcal{C}(m,n,q)I(m,n,q)
	\end{aligned}$$
	
	Note that we have
	
	$$S(N) \ll \log N \sup \limits_{C \ll Q} S(N,C)$$
	
	With a dyadic partition of $C$ fixed, we now further apply a dyadic partition on the $m$-variable. Thus we have
	
	$$S(N,C) = \sum \limits_{\substack{M_1 \ll M_0 \\ M_1 \text{ dyadic}}} S(N,C,M_1)$$

	where $S(N,C,M_1)$ is defined identically to $S(N,C)$ with the condition $m \ll M_0$ replaced by $m \sim M_1$. Note that
	
	$$S(N,C) \ll \log N \sup \limits_{M_1 \ll M_0} S(N,C,M_1)$$
	
	Therefore the bounds for $S(N)$ follows from bounds on  $S(N,C,M_1)$, which we now seek to find.
	
	\section{Conclusion of the Proof}
	
	We now evaluate the character sum and the integral getting
	
	$$\begin{aligned}
		S(N,C,M_1) = &\frac{1}{\sqrt{t}}\int \limits_{|u|\ll C/QK} \int \limits_{x\in I(M_1)} \sum \limits_{q_0^2q_1q^* = q\sim C} \frac{g(q,x)}{(q_0^2q_1q^*)^{3/2+it}}
		\\ & \times \sum \limits_{n \sim N_0}\sum \limits_{m \sim M_1} \lambda_f(m)W_{k-1}\left( \frac{2N\sqrt{m}}{q\sqrt{p}}\right) \sqrt{q^*}S(4m-pn^2,0,q_0^2)\left( \frac{\frac{4m-pn^2}{q_1}}{q^*}\right)q_1
	\end{aligned}$$
	
	Writing a general Ramanujan sum as
	
	$$S(r,0;q) = \sum \limits_{d|(r,q)}d\mu\left( \frac{q}{d}\right)$$
	
	and interchanging sums gives
	
	\begin{equation}
		\begin{aligned}
			& \frac{1}{\sqrt{t}}\int \limits_{|u|\ll C/QK} \int \limits_{x \in I(M_1)} \sum \limits_{l \ll C} \frac{\mu(l)}{l^{5/2}}\sum \limits_{dq_1q^* =q \sim C/l} \frac{g(ql,x)}{q^{5/2}}
			\\&  \times \sum \limits_{n \ll N_0} \sum \limits_{m \ll M_0} \lambda_f(m)W_{k-1}\left( \frac{4\pi N\sqrt{m}}{ql\sqrt{p}}\right) \left( \frac{d}{q^*}\right)\left( \frac{\frac{4m - pn^2}{q_1d}}{q^*}\right)q_1d\sqrt{q^*}
		\end{aligned}
	\end{equation}
	
	Applying bounds for $W_{k-1}$ and Cauchy-Schwarz we obtain
	
	$$\frac{\mu(I(M_1))}{QK\sqrt{t}}\left(1 + \frac{N\sqrt{M_1}}{C\sqrt{p}}\right)^{-1/2}\sum \limits_{l \ll C} \frac{1}{l} \left( \sum \limits_{dq_1q^* \sim C/l}\left| \sum \limits_{n\ll N_0}\sum \limits_{m \ll M_1} \lambda_f(m)\left( \frac{\frac{4m - pn^2}{q_1d}}{q^*}\right)\sqrt{q_1d}\right|^2 \right)^{1/2}$$
	
	Observe that
	
	\begin{equation}
		\mu(I(M_1))\left(1 + \frac{N\sqrt{M}}{C\sqrt{p}} \right)^{-1/2} \ll \begin{cases} CQ/N^{2-\epsilon} & \text{if }M_1 \ll C^2p/N^{2-\epsilon}
			\\ \frac{M_1^{1/4}C^{1/2}p^{1/4}}{N^{1/2}\sqrt{pK}} & \text{if } M_1 \gg C^2p/N^{2-\epsilon}
		\end{cases}
	\end{equation}

	We now apply Theorem \ref{qsieve} for each fixed $d$ and $q_1$ getting
	
	$$\frac{\mu(M_1)}{QK\sqrt{t}M_1^{1/4}} \sum \limits_{l \ll C} \frac{1}{l} \left( \sum \limits_{dq_1 \ll C/l} \left(\left(C + L\right)||\alpha(q_1,d)||_2^2 \right) \right)^{1/2}$$
	
	where
	
	\begin{equation}L \ll \sup \left\{  pN_0^2, M_1\right\}\end{equation}
	
	and 
	
	\begin{equation}
		||\alpha(q_1,d)||_2^2 = \underset{ \substack{(pn_1^2-4m_1)(pn_2^2-4m_2) = \square\\ q_1d|(pn_i^2 - 4m_i)}}{\sum \sum \sum \sum}|\lambda_f(m_1)\lambda_f(m_2)|
	\end{equation}
	
	We now apply the Deligne/Ramanujan bound to trivially bound $|\lambda_f(m_i)| \ll (pt)^\epsilon$.
	
	Now we count the number of points in the set
	
	$$\{ (n_1,n_2,m_1,m_2): n_i\sim N_0,m_i \sim M_1,(pn_1^2-4m_1)(pn_2^2-4m_2)=\square, q_1d|(pn_i^2-4m_i)\}$$
	
	Note that this is bounded by counting the set
	
	$$\{(n_1,n_2,m_1,m_2,\theta,r_1,r_2) : n_i \sim N_0,m_i \sim M_1, pn_i^2-4m_i = q_1d\theta r_i^2\}$$
	
	We observe that the total number of points $(\theta,r_1,r_2)$ is bounded by $\frac{L}{q_1d}$.
	
	Now we split into two cases - $C \ll \frac{N^{1+\epsilon}K^{1/2}}{t}$ and $C \gg \frac{N^{1+\epsilon}K^{1/2}}{t}$.
	
	\subsection{Case I - Small $q \sim C$}.
	
	Assume
	
	\begin{equation}
		C \ll \frac{N^{1+\epsilon}K^{1/2}}{t}
	\end{equation}
	
	Therefore $N_0 \ll K^{1/2}$. Trivially counting $n_1,n_2,\theta,r_2,r_2$ we obtain
	
	$$S(N,C,M_1) \ll \frac{\mu(I(M_1))}{QKt^{1/2}}\left(1 + \frac{N\sqrt{M_1}}{C\sqrt{p}} \right)^{-1/2}(C+pK)^{1/2}p^{1/2}K$$
	
	Bounding $C$ by $NK^{1/2}/t$.

	\subsection{Case II - Large $q \sim C$}
	
	We now assume $C \gg \frac{NK^{1/2}}{t}$. This implies $N_0 \sim Ct/N$ and that $n \sim N_0$ by stationary phase analysis. It also implies $pN_0^2 \gg M_0$.
	
	Here we use the crucial fact that $n \sim N_0$. We  observe that
	
	$$|pn_i^2 - \theta r_i^2| \ll M_1 \implies \left|n_i - \frac{\sqrt{\theta}r_i}{\sqrt{p}} \right| \ll \left(\frac{M_1}{pN_0} + 1 \right)$$
	
	This gives a bound on the number of $n_i$ if we sum over $r_1,r_2$ and $\theta$. So summing over $r_1,r_2$ and $\theta$ and $q_1$ and $d$ gives the point count
	
	$$\left(C + L \right)L\left( \frac{M_1}{pN_0}+1\right)^2$$
	
	If $M_1 > pN_0$ then we obtain the bound
	
	\begin{equation}
		S(N,C,M_1) \ll \frac{pN_0^2M^{3/4}}{pN_0QK\sqrt{t}} \ll \frac{N_0 M^{3/4}}{QK\sqrt{t}}
	\end{equation}
	
	Using the trivial bounds $M_1 \ll pK$ and $N_0 \ll t/\sqrt{K}$ gives
	
	\begin{equation}{\label{final1}}
		S(N,C,M_1) \ll \frac{1}{QK}p^{3/4}\sqrt{t}K^{1/4}
	\end{equation}
	
	If $M_1 \ll pN_0$ then we have the bound
	
	\begin{equation}
		S(N,C,M_1) \ll \frac{1}{QK\sqrt{t}}\mu(I(M_1))\left(1 + \frac{N\sqrt{M}}{C\sqrt{p}}\right)^{-1/2}pN_0^2
	\end{equation}
	
	Trivially bounding $N_0 \ll t/\sqrt{K}$, $M_1 \ll pK$ and $C \ll Q$ gives
	
	\begin{equation}{\label{final2}}
		S(N,M_1,C) \ll \frac{1}{QK}\frac{p^{3/4}t^{3/2}}{K^{5/4}}
	\end{equation}
	
	From $(\ref{final1})$ and $(\ref{final2})$ we get
	
	\begin{equation}{\label{finaldyadic}}
		S(N,C) \ll N^{1/2}p^{1/2}t^{1/2} \left( K^{1/4} + \frac{t}{K^{5/4}}\right)
	\end{equation}
	
	Inserting $(\ref{finaldyadic})$ into $(\ref{qdyadic})$ we obtain
	
	\begin{equation}{\label{finalbound}}
		S(N) \ll N^{1/2} p^{1/2+\epsilon}t^{1/2+\epsilon} \left( K^{1/4} + \frac{t}{K^{5/4}}\right)
	\end{equation}
	
	This proves Theorem \ref{thm2} The bound is optimal with the choice $K = t^{2/3}$ which proves Theorem \ref{thm1}.
 
\printbibliography
	
\end{document}